\documentclass{article}
\author{Dan Betea\thanks{KU Leuven, Leuven, Belgium, \texttt{dan.betea@gmail.com}. The author is supported by FWO Flanders project EOS 30889451.}}
\title{Determinantal point processes from symplectic and orthogonal characters and applications}

\usepackage{graphicx, graphics, amsmath, amsthm, amssymb, xcolor,mathtools, bbm, mathabx, hyperref, stmaryrd}
\usepackage[margin=1.2in]{geometry}

\newtheorem{thm}{Theorem}

\theoremstyle{definition}

\theoremstyle{remark}

\numberwithin{equation}{section}
\DeclareMathAlphabet{\mathpzc}{OT1}{pzc}{m}{it}

\let\P\undefined
\newcommand{\P}{\mathbb{P}}
\newcommand{\Z}{\mathbb{Z}}

\newcommand{\F}{\mathcal{F}}
\newcommand{\Gami}{\Gamma_{-}}
\newcommand{\Gapl}{\Gamma_{+}}

\newcommand{\Gsp}{\Gamma_{sp+}}
\newcommand{\Gspm}{\Gamma_{sp\pm}}
\newcommand{\Gsm}{\Gamma_{sp-}}

\newcommand{\bra}[1]{\langle #1 |}
\newcommand{\ket}[1]{| #1 \rangle}
\newcommand{\vv}{| 0 \rangle}
\newcommand{\vcv}{\langle 0 |}

\newcommand{\im }{\mathrm{i}}  

\newcommand{\dx }{\mathrm{d}}

\newcommand{\A}{\mathcal{A}_{2}}
\newcommand{\Ap}{\mathcal{A}^{+}_{2 \to 1}}
\newcommand{\Am}{\mathcal{A}^{-}_{2 \to 1}}
\newcommand{\Apm}{\mathcal{A}^{\pm}_{2 \to 1}}

\allowdisplaybreaks

\begin{document}

\maketitle

\begin{abstract}
    We show that the symplectic and orthogonal character analogues of Okounkov's Schur measure (on integer partitions) are determinantal, with explicit correlation kernels. We apply this to prove certain Borodin--Okounkov--Gessel-type results concerning Toeplitz+Hankel and Fredholm determinants; a Szeg\H{o}-type limit theorem; an edge Baik--Deift--Johansson-type asymptotical result for certain symplectic and orthogonal analogues of the poissonized Plancherel measure; and a similar result for actual poissonized Plancherel measures supported on ``almost symmetric'' partitions.
\end{abstract}

\section{Introduction and main results.}

\paragraph{Background.} Okounkov's Schur measure~\cite{oko}, along with the time-extended Schur process of Okounkov--Reshetikhin~\cite{or}, are prototypical examples of discrete determinantal point processes that have figured prominently in probability and statistical mechanics (the Kardar--Parisi--Zhang (KPZ)~\cite{Co12} universality class), combinatorics, enumerative geometry, integrable hierarchies, and other areas. The Schur measure has been extended to free boundary (the \emph{pfaffian}) cases in~\cite{br1, BBNV18}, cylindric boundaries~\cite{B2007cyl}, and also in other directions. For the combinatorial and probabilistic sides of the issue, mostly of concern to us, see the review~\cite{bg}. 

Here we are concerned with a so-called extension of the Schur measure to ``root system types other than A''. To wit, the classical Cauchy identities (orthogonality relations) for the irreducible characters of the symplectic and orthogonal groups, pairing one such function with a Schur polynomial, read as follows---see~\cite{kt,sun2} and references therein:
\begin{equation} \label{eq:cauchy_intro}
    \begin{split}
    \sum_{\lambda} sp_{\lambda} (X) s_{\lambda} (Y) &= \prod_{1 \leq i < j \leq N} (1 - y_i y_j) \prod_{1 \leq i, j \leq N} \frac{1}{(1-y_i x_j) (1-y_i x_j^{-1})}, \\
    \sum_{\lambda} o_{\lambda} (X) s_{\lambda} (Y) &= \prod_{1 \leq i \leq j \leq N} (1 - y_i y_j) \prod_{1 \leq i, j \leq N} \frac{1}{(1-y_i x_j) (1-y_i x_j^{-1})}.
    \end{split}
\end{equation}
Here $X = (x_1^{\pm 1}, \dots, x_N^{\pm 1})$, $Y = (y_1, \dots, y_N)$ are sets of variables/parameters (combinatorially, they are \emph{alphabets}). The orthogonal identity above is the one for even orthogonal groups; for the odd case one replaces $X$ by $\tilde{X} = (x_1^{\pm 1}, \dots, x_N^{\pm 1}, 1)$ and modifies the right-hand side appropriately by multiplying by $\prod_{1 \leq i \leq N} (1-y_i)^{-1}$. We can build probability measures $m_{sp/o}$ ($/$ stands throughout for \emph{or}) on integer partitions out of them via
\begin{equation} \label{eq:measures_intro}
    m_{sp}(\lambda) \ \propto \ sp_{\lambda} (X) s_{\lambda} (Y), \qquad m_{o}(\lambda) \ \propto\ o_{\lambda} (X) s_{\lambda} (Y)
\end{equation}
(viewing the $x$'s and $y$'s as numbers). These are normalized in the obvious way via~\eqref{eq:cauchy_intro}.

In this note we show that $m_{sp}$ and $m_{o}$ are determinantal measures: their $n$-point correlation functions are given by determinants of simple kernels which we write explicitly. These measures should be viewed analogously to Okounkov's Schur measure~\cite{oko}. More precisely, they are certain ``type B, C and D'' analogues of Schur measures. 

We further consider three applications: one to the study of certain Toeplitz+Hankel determinants previously encountered in random matrix theory; the second and third are asymptotic results where Kardar--Parisi--Zhang (KPZ)-type kernels and distributions appear in the edge limit: the Airy$_{2 \to 1}$ and a certain dual, as well as the Airy kernel and the classical Tracy--Widom GUE distribution.

\paragraph{Outline.} In the rest of the Introduction we present the main results. We sketch the proofs in Section~\ref{sec:proof}, and make some concluding remarks in Section~\ref{sec:rem}.

\paragraph{Main results.} We now list the main results. All but Theorem~\ref{thm:edge_schur} appear with detailed proofs in~\cite{b_sympa}. Theorem~\ref{thm:edge_schur} is new. We use throughout the standard notation language associated to integer partitions and symmetric functions as can be found in~\cite{mac}. Letters like $\lambda, \mu, \alpha, \beta$ are all reserved for partitions.

We find it convenient to use the language of specializations. A \emph{specialization} $\rho$ is just a sequence of numbers $(p_n(\rho))_{n \geq 1}$---its values on the powersum symmetric functions---assembled into the generating series $H(\rho; z) := \exp (\sum_{n \geq 1} p_{n}(\rho) z^n/n)$. The values on the complete homogeneous symmetric functions can be read from the equation $H(\rho; z) = \sum_{k \geq 0} h_{k} (\rho) z^k$. In this language, the Schur functions, (universal) symplectic, and orthogonal characters are defined via the Jacobi--Trudi $h$-formulae\footnote{There are similar formulas in terms of the elementary symmetric functions $e_n$'s. See, e.g.\,\cite{bak}. We also note we view such characters as lifted from being Laurent polynomials to the algebra of symmetric functions.}:
\begin{equation}
    \label{eq:s_jt}
    \begin{split}
    s_{\lambda}(\rho) &= \det [h_{\lambda_i-i+j}(\rho)], \\
        sp_{\lambda}(\rho) &= \frac{1}{2} \det [h_{\lambda_i - i + j}(\rho) + h_{\lambda_i - i - j + 2}(\rho)], \quad o_{\lambda}(\rho) = \det [h_{\lambda_i - i + j}(\rho) - h_{\lambda_i - i - j}(\rho)]
    \end{split}  
\end{equation}
with all determinants of size $\ell(\lambda)$ (the length of $\lambda$).

Fix two specializations $\rho^+, \rho^-$ and consider the following (possibly complex) probability measures on partitions
\begin{equation}
    m_{sp}(\lambda) = Z_{sp}^{-1} sp_{\lambda} (\rho^+) s_{\lambda} (\rho^-), \quad m_{o}(\lambda) = Z_{o}^{-1} o_{\lambda} (\rho^+) s_{\lambda} (\rho^-)
\end{equation}
where the normalization constants are the lifted versions of the right-hand side of~\eqref{eq:cauchy_intro}: $Z_{sp/o} = \exp \sum_{k \geq 1} \left( \frac{p_k (\rho^+) p_k (\rho^-)}{k} \pm \frac{p_{2k} (\rho^-)}{2k} - \frac{p_k^2(\rho^-)}{2k} \right)$ (we take $+$ for $sp$ and $-$ for $o$). Our main result is as follows.

\begin{thm} \label{thm:corr}
    \label{thm:corr_lifted}
    Assume that the specializations $\rho^{\pm}$ are such that $k^{-1} p_k(\rho^{\pm}) = O(r_{\pm}^k)$ for some constants $r_\pm > 0$ with $\min(1, 1/r_+) > r_- > 0$. Fix $k_i \in \Z+1/2\,\forall 1 \leq i \leq n$. Then the measures $m_{sp}$ and $m_o$ are determinantal:
    \begin{equation}
        m_{sp/o} (\{k_1, \dots, k_n\} \subset \{ \lambda_i - i + 1/2\}) = \det [K_{sp/o} (k_i, k_j)]_{1 \leq i,j \leq n}
    \end{equation}
    with correlation kernels given by 
\begin{equation}
    \label{eq:kernel_lifted}
    \begin{split}
    K_{sp}(a, b) &= \int_{z} \int_{w} \frac{F(z)}{F(w)} \frac{(1-w^2)}{(1-wz) (1-wz^{-1})} \frac{\dx z \dx w}{(2 \pi \im)^2 z^{a+3/2} w^{-b+1/2}}, \\
    K_{o}(a, b) &= \int_{z} \int_{w} \frac{F(z)}{F(w)} \frac{(1-z^2)}{(1-wz) (1-wz^{-1})} \frac{\dx z \dx w}{(2 \pi \im)^2 z^{a+3/2} w^{-b+1/2}}
    \end{split}
\end{equation}
where $F(z) := \frac{H(\rho^+; z)}{H(\rho^-; z) H(\rho^-; z^{-1})}$
and the contours are simple closed counterclockwise curves around 0 such that $\min(1/|w|, 1/r_+, 1/r_-) > |z| > |w| > r_- > 0$.
\end{thm}

The contour integrations above, while analytical, have the simple combinatorial meaning of coefficient extraction. In other words what we really obtain are the generating series for the kernels. Nonetheless, they are quite useful in performing analytical computations (like steepest descent analysis). 

The first application we consider is to the study of certain Toeplitz+Hankel determinants which appeared in random matrix theory~\cite{joh5} (and which are similar to the ones of Baik--Rains~\cite{br1} who study longest increasing subsequences of random involutions). 

Start with a function $f(z) = \exp(R_+(z) + R_-(z)) = \sum_{k \in \Z} f_k z^k$ (the \emph{symbol})\footnote{The $f_k$'s are the Fourier coefficients of $f$, and $f = \exp R_+ \exp R_-$ is its \emph{Wiener--Hopf} factorization.}, where $R_{\pm}(z) = \sum_{k \geq 1} \rho_k^{\pm} z^{\pm k}$ for appropriately chosen numbers $(\rho^{\pm}_k)_{k \geq 1}$ (chosen so that the exponents in the right-hand sides of Theorem~\ref{thm:szego} are absolutely convergent). Let $\check{f}(z) := f(-z)^{-1}$ having Fourier coefficients $(\check{f}_k)_{k \in \Z}$. Consider the following Toeplitz+Hankel determinants:
\begin{equation}
    \begin{split}
        D^1_{n} = \det [f_{i-j} + f_{i+j}]_{0 \leq i,j \leq n-1}, \qquad D^2_{m} = \det [\check{f}_{i-j} - \check{f}_{i+j+2}]_{0 \leq i,j \leq m-1}, \\ 
        D^3_{n} = \det [f_{i-j} - f_{i+j+2}]_{0 \leq i,j \leq n-1}, \qquad D^4_{m} = \det [\check{f}_{i-j} + \check{f}_{i+j}]_{0 \leq i,j \leq m-1}.
    \end{split}
\end{equation}  

These determinants can be written as observables for the measures $m_{sp/o}$, and lead to the following symplectic/orthogonal variants of the Gessel formulae of~\cite{ges}.

\begin{thm} \label{thm:gessel}
It holds that:
\begin{equation}
    \begin{split}
        \frac{1}{2} D^1_n = \sum_{\lambda: \ell(\lambda) \leq n} sp_{\lambda} (\rho^+) s_{\lambda} (\rho^-), \qquad D^2_m = \sum_{\lambda: \lambda_1 \leq m} sp_{\lambda} (\rho^+) s_{\lambda} (\rho^-), \\
        D^3_n = \sum_{\lambda: \ell(\lambda) \leq n} o_{\lambda} (\rho^+) s_{\lambda} (\rho^-), \qquad  \frac{1}{2} D^4_m = \sum_{\lambda: \lambda_1 \leq m} o_{\lambda} (\rho^+) s_{\lambda} (\rho^-).
        \end{split}
\end{equation}
\end{thm}

Furthermore, a consequence of the above and the (lifted) Cauchy identities is the following asymptotic Szeg\H{o}-type result.

\begin{thm} \label{thm:szego}
It holds that:
\begin{equation}
    \begin{split}
        \lim_{n \to \infty} \frac{1}{2} D^1_n = \lim_{m \to \infty} D^2_m &= \exp \sum_{k \geq 1} \left( k \rho^+_k \rho^-_k + \rho^-_{2k} - \frac{k(\rho^-_k)^2}{2} \right), \\
        \lim_{n \to \infty} D^3_n = \lim_{m \to \infty}  \frac{1}{2}  D^4_m &= \exp \sum_{k \geq 1} \left( k \rho^+_k \rho^-_k - \rho^-_{2k} - \frac{k(\rho^-_k)^2}{2} \right).
    \end{split}
\end{equation}
\end{thm}

Using different methods, the result was first proved by Johansson~\cite{joh5}, and reproved subsequently by Basor--Ehrhardt, Dehaye and Deift--Its--Krasovsky in~\cite{be4, deh, dik}. 

We further have a Borodin--Okounkov-type~\cite{bo} result of the form Toeplitz+Hankel = Fredholm---see also~\cite{tw2}. It appears in~\cite{be4} in a different purely analytical form. We note that the right-hand sides below are \emph{Fredholm determinants}.
\begin{thm} \label{thm:bo}
        It holds that:
        \begin{equation}
            D_m^2 = Z_{sp} \cdot \det (1-K_{sp})_{\ell^2 \{ m+\frac{1}{2}, m+\frac{3}{2}, \dots \} }, \quad \frac{1}{2} D_m^4 = Z_{o} \cdot \det (1-K_o)_{\ell^2 \{ m+\frac{1}{2}, m+\frac{3}{2}, \dots \} }
        \end{equation}
        where $Z_{sp/o}$ are as above; and where $K_{sp/o}$ are as in Theorem~\ref{thm:corr} and we assume the same conditions on $\rho^{\pm}$ as in thm.\,cit.
\end{thm}

Similar results hold for $D^1$ and $D^3$ via duality between largest parts and lengths of partitions. We omit the details. 

We remark that all of the above is purely algebraic combinatorial. I.e.,\,if one removes all the analytical convergence conditions, the statements still make sense at the level of formal power series. In particular it is entirely possible that the analytical conditions we impose for convergence can be relaxed.

As a different type of application, for a positive real $\theta$, we consider the measures 
\begin{equation}
    P_{sp} (\lambda) \ \propto \ sp_{\lambda} (pl_{2 \theta}) s_{\lambda} (pl_{\theta}), \qquad P_{o} (\lambda) \ \propto \ o_{\lambda} (pl_{2 \theta}) s_{\lambda} (pl_{\theta})
\end{equation}
where $pl_{\theta}$ is the Plancherel specialization sending the first powersum $p_1$ to $\theta$ and all the rest to $0$. These are \textit{signed} measures (in some sense, analogues of Okounkov's Schur version of the poissonized Plancherel measure on partitions) whose edge asymptotic behavior (as $\theta \to \infty$) we analyze in detail. See~\cite{b_sympa} for the bulk behavior as well. We believe that (at the edge) the behavior of these measures is universal, common to all symplectic/orthogonal Schur measures having certain properties. The following theorem is then an analogue of the Baik--Deift--Johansson~\cite{bdj} theorem on (poissonization of) the longest increasing subsequence of random permutations.

\begin{thm} \label{thm:edge}
We have:
    \begin{equation}
    \lim_{\theta \to \infty} P_{sp}\left(\frac{\lambda_1 - 2 \theta}{\theta^{1/3}} \leq s\right) = \det(1-\Ap), \quad \lim_{\theta \to \infty} P_{o}\left(\frac{\lambda_1 - 2 \theta}{\theta^{1/3}} \leq s\right) = \det(1-\Am)
    \end{equation}
    where the Fredholm determinants on the right are on $L^2(s, \infty)$ of the operators given by the following kernels:
    \begin{equation}
        \Apm(x, y) := \int_{0}^{\infty} Ai(x+s) Ai(y+s) \dx s \pm \int_{0}^{\infty} Ai(x-s) Ai(y+s) \dx s
    \end{equation}
    with $Ai$ is the Airy Ai function (solution of $y''=xy$ decaying exponentially at $\infty$).
\end{thm}

We can prove more: for $n = 1,2,3,4\dots$ fixed, the largest first $n$ points of $\{\lambda_i - i \}$ in the appropriate (edge) scaling of Theorem~\ref{thm:edge} converge, in the sense of finite dimensional distributions, to the first $n$ points of the $\Ap$ process for $P_{sp}$ and the $\Am$ process for $P_o$. We omit the statement for brevity.

Let us remark that $\Ap$ is the Airy${}_{2 \to 1}$ kernel of Borodin--Ferrari--Sasamoto~\cite{bfs}. It appears in the scaling limit of TASEP with half-flat initial conditions---see also~\cite{bz2}. The associated distribution (Fredholm determinant) above interpolates between the Tracy--Widom~\cite{tw} GUE distribution and the corresponding GOE~\cite{TW05} distribution---see~\cite{bfs}.

The kernel $\Am$ appears to be new, and could be described as the \emph{dual} (in the sense of symmetric functions) of $\Ap$.

Finally, let us consider the following sets of almost-symmetric partitions:
\begin{equation}
    A = \{ \alpha = (a_1, a_2, \dots|a_1+1, a_2+1, \dots) \}, \quad B = \{ \beta = (b_1+1, b_2+1, \dots|b_1, b_2, \dots) \}
\end{equation}
where the coordinates above are~\textit{Frobenius coordinates}~\cite{mac}. Clearly the partitions in $B$ are the transposed of those in $A$. Fix again a $\theta > 0$ and consider the (actual) probability measures on $A, B$ given by the poissonized dimensions of the respective (symmetric group) representations (i.e.\,$\dim \lambda$ is the number of standard Young tableaux of shape $\lambda$):
\begin{equation}
    \P_A(\alpha) = e^{-\theta^2/2} \frac{\theta^{|\alpha|} \dim \alpha}{|\alpha|!}, \quad \P_B(\beta) = e^{-\theta^2/2} \frac{\theta^{|\beta|} \dim \beta}{|\beta|!}.
\end{equation}
(Note the weights are identical, but the supports are different for $\P_A, \P_B$.) We have the following Baik--Deift--Johansson-type asymptotic result for the largest part of such partitions.

\begin{thm} \label{thm:edge_schur}
        We have:
            \begin{equation}
            \lim_{\theta \to \infty} \P_A \left(\frac{\alpha_1 - 2 \theta}{\theta^{1/3}} \leq s\right) = \lim_{\theta \to \infty} \P_B \left(\frac{\beta_1 - 2 \theta}{\theta^{1/3}} \leq s \right) = F_2(s) := \det(1-\A)
            \end{equation}
            where $F_2(s)$ is the Tracy--Widom GUE distribution~\cite{tw} given as the Fredholm determinant on $L^2(s, \infty)$ of the Airy kernel $\A(x, y) := \int_{0}^{\infty} Ai(x+s) Ai(y+s) \dx s$.
\end{thm}

As above, we can state and prove a more general result showing convergence of the largest (few) parts of such partitions to the corresponding particles in the Airy ensemble. We omit the statements for brevity.

\section{A sketch of the proofs} \label{sec:proof}

\paragraph{Proofs of Theorems~\ref{thm:gessel}, \ref{thm:szego} and~\ref{thm:bo}.} The argument that proves Theorem~\ref{thm:gessel} is similar to that used by Gessel~\cite{ges}. It is a simple application of the Cauchy--Binet identity together with the Jacobi--Trudi (both the $h$ and $e$ ones) formulae for the Schur, symplectic, and orthogonal functions. The missing link is to observe, as Tracy--Widom have done for the Schur measure~\cite{tw2}, that the Fourier coefficients of the logarithm of $f$, the numbers $\rho^\pm_k$ for $k \geq 1$, naturally define specializations $\rho^\pm$ via the formulas $p_k(\rho^\pm) = \rho^{\pm}_k$.

Theorem~\ref{thm:szego} is just a consequence of the symplectic and orthogonal Cauchy identities (the identities~\eqref{eq:cauchy_intro} where we lift the $X, Y$ parameters to more generic specializations by just replacing the powersums $p_n(X), p_n(Y)$ by the ones for the specializations). The \emph{restricted Cauchy sums} on the right-hand side of the equations of Theorem~\ref{thm:gessel} just go to the appropriate \emph{unrestricted} Cauchy sums in the limit, and the latter are then equal to the Cauchy kernels $Z_{sp/o}$ (by definition the right-hand sides of Theorem~\ref{thm:szego}, and nothing else than the lifts of the right-hand sides of~\eqref{eq:cauchy_intro} to arbitrary specializations).

Finally, Theorem~\ref{thm:bo} is based on the following observation of Borodin--Okounkov~\cite{bo}: on one hand, $D^2$ and $D^4$ are gap probabilities (that no parts of $\lambda$ are bigger than $m$) of the appropriate symplectic/orthogonal measures via the Gessel-like formulae of Theorem~\ref{thm:gessel}. On the other hand, since by our main Theorem~\ref{thm:corr} such measures are determinantal, gap probabilities are, by inclusion--exclusion, \emph{Fredholm determinants} (up to the normalizing constants $Z_{sp/o}$ needed to make the restricted Cauchy sums actual probabilities).

\paragraph{Proof of Theorem~\ref{thm:corr}.} We now sketch the proof of our main theorem. Aided by a construction of~\cite{bak} and the analogous Schur result and techniques of Okounkov~\cite{oko}, we use infinite wedge formalism to construct our correlation functions (observables) as certain expectations inside fermionic Fock space. A Wick lemma then gives the determinantal structure. To sketch the proof, we concentrate on the symplectic measure. The orthogonal case can be arrived at \emph{combinatorially} by applying the symmetric function involution $\omega$ sending $\omega(h_n) = e_n$, as in this case $\omega(sp_\lambda) = o_{\lambda'}$. We recall below basic facts about fermionic Fock space, see~\cite{oko, bak, b_sympa} for more details.


We denote by $\F$ the \emph{fermionic Fock space}, the Hilbert space spanned by the orthonormal basis $\ket{S}$, where $S
 \subset \Z' := \Z + 1/2$ with $|S \cap \Z'_+|<\infty$ and $|\Z'_-  - S| < \infty$. Elements of $S$ are \textit{particles}. We use  bra--ket notation throughout. A basis vector $\ket{S}$ is the semi-infinite wedge product
$\ket{S} = \underline{s_1} \wedge \underline{s_2} \wedge \underline{s_3} \wedge \cdots$
where $s_1 > s_2 > s_3 > \cdots$ are \textit{particle positions}. States $S$ are in 1--1 correspondence with pairs $(\lambda, c)$, where $\lambda$ is a partition and $c \in \Z$ is the \emph{charge}. The correspondence is $(\lambda, c) \longleftrightarrow S(\lambda,c):=\{ \lambda_i - i + 1/2 + c ,\ i \geq 1\}$. It implies the grading $\F = \oplus_{c \in \Z} \F_c$ where $\F_c$ is spanned by partitions $\lambda$ of charge $c$, which we denote by $\ket{\lambda, c}$. Importantly, we omit the $c$ if $c=0$ and write $\ket{\lambda}$ for $\ket{\lambda, 0} \in \F_0$. In fact, we are only interested in the subspace $\F_0$. The vector $\vv$, corresponding to the empty partition (of charge 0), is called the \emph{vacuum} vector. 

For $k \in \Z'$, define the free fermionic creation operator $\psi_k \ket{S} := \underline{k} \wedge \ket{S}$ with $\psi^*_k$ its adjoint under the inner product $\langle \lambda, c | \mu, d \rangle = \delta_{\lambda, \mu} \delta_{c, d}$. Roughly speaking, $\psi_k$ tries to add a particle to the configuration $S$, while $\psi^*_k$ tries to remove one (with an appropriate sign in both cases). Both kill the vector if a particle is there at site $k$ (for $\psi_k$) or is absent (for $\psi^*_k$). These operators satisfy the canonical anti-commutation relations (CAR) $\{ \psi_k, \psi_\ell^* \} = \delta_{k,\ell}, \{ \psi_k, \psi_\ell \} = \{ \psi_k^*, \psi_\ell^* \} = 0, k,\ell \in \Z'$ and we collect them into generating series $\psi(z) := \sum_{k \in \Z'} \psi_k z^k, \psi^*(w) := \sum_{k \in \Z'} \psi^*_k w^{-k}$. 
We build the \emph{bosonic operators} $\alpha_n$ out of fermions as follows: $\alpha_n := \sum_{k \in \Z'} \psi_{k-n} \psi_k^*$ for $n = \pm 1, \pm 2, \ldots$ 
We have $\alpha_n^*=\alpha_{-n}$, and $\alpha_n \vv=0, n>0$. These operators satisfy the canonical commutation relations (CCR) $[\alpha_n,\alpha_m] = n \delta_{n,-m}$ and commute thusly with the fermionic fields: $[\alpha_n,\psi(z)] = z^n \psi(z), [\alpha_n,\psi^*(w)] = - w^n \psi^*(w)$. Further define the following so-called \emph{half-vertex operators} $\Gamma_\pm(\rho), \Gspm(\rho)$ (for $\rho$ a specialization) by
\begin{equation}
  \Gamma_\pm(\rho) := \exp \left( \sum_{n \geq 1} \frac{p_n(\rho) \alpha_{\pm n}}{n} \right), \quad \Gspm(\rho) := \exp \sum_{n \geq 1} \frac{1}{n} \left( p_n(\rho) \alpha_{\pm n} + \frac{\alpha_{\pm 2n}}{2} - \frac{\alpha_{\pm n}^2}{2} \right).
\end{equation}
They have the property of having (skew) Schur functions and (lifted) symplectic characters as matrix elements (see~\cite{bak} for the latter), as follows:
\begin{equation}
    s_{\mu/\lambda}(\rho) = \bra{\lambda} \Gamma_+(\rho) \ket{\mu} = \bra{\mu} \Gamma_-(\rho) \ket{\lambda}, \quad sp_{\lambda}(\rho) = \vcv \Gsp(\rho) \ket{\lambda} = \bra{\lambda} \Gsm(\rho) \vv.
\end{equation}
and commute among themselves and with the fermionic fields $\psi, \psi^*$ thusly\footnote{We only list the relevant commutations for our purposes.}:
\begin{align}
    & \Gsp(\rho^+) \Gami(\rho^-) = H(\rho^+; \rho^-) h_{sp}(\rho^-) \Gami(Y) \Gapl^{-1} (\rho^-) \Gsp(\rho^+), \notag \\
    & \Gamma_\pm(\rho^{-}) \psi(z) = H(\rho^{-};z^{\pm 1}) \psi(z) \Gamma_\pm(\rho^{-}), \quad \Gamma_\pm(\rho^{-}) \psi^*(w) = H(\rho^{-}; w^{\pm 1})^{-1} \psi^*(w) \Gamma_\pm(\rho^{-}), \notag \\
    & \Gsp(\rho^+) \psi(z) \Gsp^{-1}(\rho^+) = H(\rho^+; z) \psi(z) \Gapl^{-1}(z), \\
    & \Gsp(\rho^+) \psi^*(w) \Gsp^{-1}(\rho^+) = (1-w^2) H(\rho^+; w)^{-1} \psi^*(w) \Gapl(w) \notag
\end{align}
where $H(\rho; \rho') := \exp \sum_{n \geq 1} p_n(\rho) p_n(\rho')/n$ and $h_{sp}(\rho) := \exp \sum_{n \geq 1} \left( \frac{p_{2n}(\rho)}{2n} - \frac{p^2_{n}(\rho)}{2n} \right)$. We note in passing the first relation above already implies the symplectic Cauchy identity.

To find the correlation functions, starting from the simple observation that $\psi_k \psi_k^* \ket{\lambda}$ picks out only the configurations $\ket{\lambda}$ satisfying $k \in \{ \lambda_i - i + 1/2 \}$, we can write the $n$-point correlation functions as $\frac{1}{Z} \left[ \frac{z_1^{k_1} \dots z_n^{k_n}}{w_1^{k_1} \dots w_n^{k_n}} \right] \left\langle \Gsp(\rho^+) \left( \prod_{i=1}^n  \psi(z_i) \psi^* (w_i) \right) \Gami(\rho^-) \right\rangle$ where $[\cdot]$ stands for coefficient extraction and we used $\langle A \rangle := \vcv A \vv$. Moreover, we have the Wick lemma\footnote{Recall the time ordering notation $\mathcal{T} (\psi(z_i) \psi^* (w_j)) = \psi(z_i) \psi^* (w_j)$ if $i \leq j$ and $= - \psi^* (w_j) \psi(z_i)$ otherwise.}
\begin{equation}
    \left\langle \Gsp(\rho^+) \left( \prod_{i=1}^n  \psi(z_i) \psi^* (w_i) \right) \Gami(\rho^-) \right\rangle = \det_{i, j} \left\langle \Gsp(\rho^+) \mathcal{T} (\psi(z_i) \psi^* (w_j)) \Gami(\rho^-) \right\rangle
\end{equation}
which follows from the determinantal evaluation 
\begin{equation}
\det_{i,j} \frac{1}{(1-w_i z_j)(1-w_i z_j^{-1})} =  \frac{\prod_{i<j} (z_i - z_j) (w_i - w_j) (1-w_i w_j) (1- z_i^{-1} z_j^{-1})} { \prod_{i,j} (1-w_i z_j) (1-w_i z_j^{-1})}.
\end{equation}
The matrix elements inside the determinant can be evaluated using commutation relations above between the $\Gamma$ operators and the $\psi(z),\psi^*(w)$'s leading to the stated generating series of our kernels. The fact that coefficient extraction can be turned into contour integrals is routine and this finishes the proof.


\paragraph{Proof of Theorem~\ref{thm:edge}.} The proof goes via classical steepest descent analysis. It consists of two parts. The first part is to show that the kernels for $P_{sp/o}$, determinantal with $F(z) = \exp \theta (z-z^{-1})$ in the notation of Theorem~\ref{thm:corr}, converge to the corresponding $\Apm$ kernels. For the edge limit, we scale the particle positions $a, b$ in the kernel $K(a, b)$ as $(a, b) = 2 \theta + \theta^{1/3} (x, y)$ as $\theta \to \infty$. In this asymptotic regime, it can be proven via a by now standard (almost ``folklore'') limiting argument that in the integral representation of $K(a, b)$ the only non-vanishing contribution comes from a neighborhood of size $\theta^{1/3}$ of the double critical point $z=1$ of the action $S(z) = z - z^{-1} - 2 \log z$. Scaling then the integration variables as $(z, w) = (e^{\zeta \theta^{-1/3}}, e^{\omega \theta^{-1/3}})$, this contribution becomes finite. It becomes the standard Airy 2-to-1 kernel of~\cite{bfs} in the case of the $sp$ measure (and $\Am$ otherwise), in a slightly different representation written as a double contour integral---see~\cite{b_sympa, bfs} for the actual formula and the first reference for more details (along with another slightly different proof starting from Bessel-like kernels). 

The second part is to show that the finite-size discrete Fredholm determinants converge to the Fredholm determinants of the Airy-like operators. This uses the same analysis as above plus a Hadamard-type bound on determinants. The argument is standard, see e.g.\,~\cite{boo}. Moreover, the same argument can be used to prove multi-point results: the first few parts of the discrete models converge to the first few parts in the continuous ensemble---see again~\cite{boo} for the classical Plancherel measure, but the structure of the argument is the same.

\paragraph{Proof of Theorem~\ref{thm:edge_schur}.} We first show how to obtain said Plancherel measures from symplectic and orthogonal ones. First one has (see e.g.\,\cite{bak} and references therein)
\begin{equation} \label{eq:sp_o_s}
    sp_{\lambda}(\rho) = \sum_{\alpha \in A} (-1)^{|\alpha|/2} s_{\lambda / \alpha}(\rho), \quad o_{\lambda}(\rho) = \sum_{\beta \in B} (-1)^{|\beta|/2} s_{\lambda / \beta}(\rho).
\end{equation}
Putting $\rho = 0$ gives $sp_\lambda = (-1)^{|\lambda|/2}$ if $\lambda \in A$ and $0$ otherwise (similarly for $o$ and $B$). Thus we choose specializations for $m_{sp/o}$ as $\rho^+ = 0, \rho^- = pl_{i \theta}$ (the imaginary unit is needed to cancel the sign just above). The measures are then determinantal with correlations given by Theorem~\ref{thm:corr} with $F(z) = \exp[-i \theta(z+z^{-1})]$. Changing $-i z \mapsto z$ (and similarly for $w$) changes the action to the one from the previous proof and then all the above arguments apply. One difference is that under this latter rotation, the factor (say for the $sp$ case) $\frac{1-w^2}{(1-wz) (1-w/z)}$ inside the integral changes to $\frac{1+w^2}{(1+wz) (1-w/z)}$ and then only the term $(1-w/z)^{-1}$ contributes asymptotically. This is responsible for the ``collapse'' to the Airy 2 kernel (as opposed to the more general ones from Theorem~\ref{thm:edge}).

\section{Concluding remarks} \label{sec:rem}

In this paper we showed how symplectic and orthogonal characters lead to (a priori complex) probability measures which are in turn determinantal point processes. We have also shown how these results have interesting finite and asymptotic applications. We make a few concluding remarks.

First, other than the asymptotic statements of Theorems~\ref{thm:edge} and~\ref{thm:edge_schur}, the other proofs are or can be made combinatorial. Furthermore, the Szeg\H{o} limit Theorem~\ref{thm:szego} can even be made bijective. For specializations \emph{into variables} and given the Gessel formula (itself combinatorially following from Cauchy--Binet), it is equivalent to the classical Cauchy identity from~\eqref{eq:cauchy_intro}. The latter has a purely combinatorial proof in terms of tableaux insertion algorithms. For the symplectic version of such algorithms see Sundaram's thesis~\cite{sun1}.

Second, consider the following Cauchy--Littlewood identities~\cite{mac}, proven bijectively (using insertion algorithms) by Burge~\cite{bur}:
\begin{equation}
    \sum_{\lambda \in A} s_\lambda (x_1, \dots, x_N) = \prod_{i < j} (1+x_i x_j), \quad \sum_{\lambda \in B} s_\lambda (x_1, \dots, x_N) = \prod_{i \leq j} (1+x_i x_j).
\end{equation}
They lead to natural measures on partitions, and such measures are determinantal (first proved by Rains~\cite[Section 7]{rai}, the result is also a corollary of Theorem~\ref{thm:corr}). The case $x_i = \theta/N, N \to \infty$ leads to our Plancherel measures $\P_{A/B}$. More mysteriously, the case of all $x_i = 1$ (in fact even the generic case) yields a bijection between size $N$ or $N+1$ Aztec diamonds and certain Young tableaux (symmetric plane partitions with diagonal slice in $A$ or $B$), as \emph{both} sets of objects are in bijection with $0-1$ triangles of numbers: one bijection is via the Burge correspondences~\cite{bur}, the other via a growth diagram version of the shuffling algorithm as described in~\cite{bbbccv}. It would be interesting to see how one can exploit this. The asymptotic analysis of the general such measures would also be interesting to perform.

Third, the measures $\P_{A/B}$ can be viewed, using the terminology of~\cite{BBNV18}, as one free-boundary Schur measures, and their asymptotic analysis seems to complete the one for such measures introduced in the work of Baik--Rains~\cite{br1, br2}. See also~\cite{BBNV18} for a similar take on the rest of such free boundary measures. 

Fourth and connected to the previous paragraph, it would be interesting to give measures $\P_{A/B}$ (along with the generalizations above-described) corner-growth/longest increasing subsequence interpretations similar to the classical ones of Baik--Rains~\cite{br1}. In fact, this would be interesting in full generality for symplectic and orthogonal measures. One partial answer to the latter can be found in~\cite{bz1}. What seem to be missing are appropriate Greene (or even Schensted) theorems~\cite{gre}.

\paragraph{Acknowledgements.} The author is grateful to Estelle Basor, J\'er\'emie Bouttier, Elia Bisi, Alexei Borodin, Sasha Bufetov, L\'aszl\'o Erd\H{o}s, Patrik Ferrari, Christian Krattenthaler, Peter Nejjar, Anita Ponsaing, Eric Rains, Arun Ram, Craig Tracy, Mirjana Vuleti\'c, Michael Wheeler, Paul Zinn-Justin and Nikos Zygouras for illuminating conversations and references to the literature. Part of this work was performed while the author was visiting the University of Melbourne, whose hospitality is gratefully acknowledged.

\bibliographystyle{myhalpha}
\bibliography{symplectique_FPSAC2020}

\end{document}